\documentclass[12pt, a4paper]{article}
\usepackage[top=2.5cm,bottom=2.5cm,left=2.5cm,right=2.5cm]{geometry}
\usepackage[russian,english]{babel}
\usepackage{amsmath}
\usepackage[standard,thmmarks,amsmath]{ntheorem}
\usepackage{amssymb}
\usepackage{graphicx}
\usepackage[normalem]{ulem}

\newtheorem{note}{Note}
\def\le{\leqslant}

\def\ge{\geqslant}

\title{\vspace*{2.5cm}\bf Minimizing Wiener Index for Vertex-Weighted Trees with Given Weight and Degree Sequences}
\author{\bf Mikhail Goubko\thanks{This research is supported by the grant of Russian
Foundation for Basic Research, project No 13-07-00389.}\\
\emph{Institute of Control Sciences}\\
\emph{of Russian Academy of Sciences, Moscow, Russia}\\
\texttt{mgoubko@mail.ru}}
\date{\vspace{0.1in} \normalsize (Received Aug 20, 2014)}

\begin{document}

\maketitle

\thispagestyle{empty}

\begin{abstract}
In 1997 Klav\v{z}ar and Gutman suggested a generalization of the Wiener index to vertex-weighted graphs. We minimize
the Wiener index over the set of trees with the given vertex weights' and degrees' sequences and show an optimal tree
to be the, so-called, Huffman tree built in a bottom-up manner by sequentially connecting vertices of the least
weights.
\end{abstract}

\baselineskip=0.30in

\section{Introduction}

In 1947 Harold Wiener \cite{Wiener47} employed the sum of distances between vertices in a chemical graph representing a
molecule to explain boiling points of alkanes. Later the sum of distances between all vertices in a graph was called
the Wiener index, which became one of the earliest topological indices.

Since then extensive research was performed on revealing connection between different topological indices of molecules
and physical, chemical, pharmacological, and biological properties of substances (see, for instance, \cite{Balaban}),
and the Wiener index appeared to be among the most useful and powerful ones (see \cite{Dobrynin01}).

For a simple connected undirected graph $G$ with the vertex set
$V(G)$ and the edge set $E(G)$ and for any pair of vertices $u,v\in
V(G)$ let $d_G(u, v)$ denote the distance (the length of the
shortest path) between $u$ and $v$ in $G$. Then the Wiener index of
the graph $G$ is defined as
$$WI(G):=\frac{1}{2}\sum_{u,v\in V(G)}d_G(u,v).$$

In 1997 Klav\v{z}ar and Gutman \cite{Klavzar97} suggested a
generalization of the Wiener index to vertex-weighted graphs. They
endowed each vertex $v\in V(G)$ in graph $G$ with some weight
$\mu_G(v)$ (in contrast to integer weights, originally used in
\cite{Klavzar97}, below we allow for arbitrary non-negative weights)
and defined the \emph{vertex-weighted Wiener index} (\emph{VWWI})
for such a graph as

$$VWWI(G):=\frac{1}{2}\sum_{u,v\in V(G)}\mu_G(u)\mu_G(v) d_G(u,v).$$

When the weight of each vertex in a graph $G$ is equal to the degree of this vertex in $G$, this index is referred to
as the \emph{Schultz index of the second kind} \cite{Gutman94} or the \emph{Gutman index} \cite{Todeschini00}.

One of the typical problems in topological index study is estimation of index value bounds over the certain class of
graphs (molecules). In \cite{Volkmann02} a tree, which minimizes the Wiener index over the set of all trees with the
given maximum vertex degree $\Delta$ has been shown to be a balanced $\Delta$-tree (the, so-called, \emph{Volkmann
tree}). Lin \cite{Lin13}, and Furtula, Gutman, and Lin \cite{Furtula13} explored minimizers and maximizers of the
Wiener index for trees of the fixed order and all degrees odd. Wang \cite{Wang08} and Zhang et al. \cite{Zhang08} have
shown independently that the minimizer of the Wiener index over the set of trees with the given vertex degrees'
sequence is the, so-called, \emph{greedy tree} \cite{Wang08}. It is built in top-down manner by adding vertices from
the highest to the lowest degree to the seed (a vertex of maximum degree) to keep the tree as balanced as possible.

In the present paper we extend the results of \cite{Wang08, Zhang08} to the vertex-weighted trees and show that some
generalization of the famous Huffman algorithm \cite{Huffman52} for the optimal prefix code builds an optimal tree,
which coincides with the greedy tree in  case of unit weights.

The paper has the following structure. In Section \ref{section_Huffman_alg} we describe the generalized Huffman
algorithm and announce the main theorem. In Section \ref{section_huffman_tree} we immerse the problem of index
minimization into the space of directed trees, which is more convenient to study. We define the notion of the vector of
subordinate groups' weights playing the key role in the proofs, and prove some important properties of Huffman trees.
In Section \ref{section_majorization} we follow the line of the proof from \cite{Zhang08} establishing the relation
between index minimization and the majorization problem of vectors of subordinate groups' weights. In Section
\ref{section_minimizing_index} we introduce the notion of a proper tree and combine the above results proving that the
Huffman tree minimizes \emph{VWWI}. We discuss possible extensions in the concluding section.

\section{Wiener Index and Huffman Trees}\label{section_Huffman_alg}
\subsection{Generating Tuples}
For a simple connected undirected graph $G$ and a vertex $v\in V(G)$ let us denote with $d_G(v)$ its \emph{degree},
i.e., the number of vertices being incident to $v$ in $G$. Denote with $W(G)$ the set of \emph{pendent vertices} (those
having degree 1) of the graph $G$, and with $M(G):=V(G)\backslash W(G)$ the set of \emph{internal vertices} (with
degree greater than unity) of $G$.

\begin{definition}
A simple undirected graph $G$ is called \emph{vertex-weighted} if
each vertex $v\in V(G)$ is endowed with a non-negative number
$\mu_G(v)>0$. The total vertex weight of the graph $G$ is denoted
with $\mu_G$.
\end{definition}

A connected vertex-weighted graph $T$ with $N$ vertices and $N-1$ edges is called a \emph{vertex-weighted tree}. Denote
with $\mathcal{T}$ the set of all vertex-weighted trees.

All graphs below are supposed to be vertex-weighted, unless stated otherwise.

\begin{definition}
Consider a vertex set $V$. Let the function $\mu: V \rightarrow
\mathbb{R}_+$ assign a non-negative weight $\mu(v)$ to each vertex
$v\in V$, while the function $d: V \rightarrow \mathbb{N}$ assigning
a natural degree $d(v)$. The tuple $\langle\mu, d\rangle$ is called
\emph{a generating tuple} if the following identity holds:
\begin{equation}\label{tree_degrees_identity}
\sum_{v\in V} d(v)=2(|V|-1).
\end{equation}

Let $\mathcal{T}(\mu, d):=\{T\in \mathcal{T}: V(T)=V, d_T(v)=d(v), \mu_T(v)=\mu(v) \text{ for all }v\in V\}$ be the set
of trees with the vertex set $V$ and vertices having weights $\mu(v)$ and degrees $d(v)$, $v\in V$. Also denote with
$\overline{\mu}:=\sum_{v\in V}\mu(v)$ the total weight of the vertex set $V$.
\end{definition}

It is well-known that $\mathcal{T}(\mu, d)$ is not empty if and only if identity (\ref{tree_degrees_identity}) holds.

Let $V(\mu, d)$ be the domain of functions of a generating tuple $\langle\mu, d\rangle$. Introduce the set $W(\mu,
d):=\{w \in V(\mu, d): d(w)=1\}$ of \emph{pendent} vertices and the set $M(\mu, d):=V(\mu,d)\backslash W(\mu,d)$ of
\emph{internal} vertices.

Below we refer to the typical generating tuple as $\langle\mu, d\rangle$, which is defined on the vertex set $V:=V(\mu,
d)$ with the pendent vertex set $W:=W(\mu, d)$ consisting of $n=|W|\ge 2$ vertices and the internal vertex set
$M:=M(\mu, d)$ consisting of $q=|M|\ge 1$ vertices.

We will solve the problem of characterizing the set
$$\mathcal{T}^*(\mu, d):=\text{Argmin}_{T\in
\mathcal{T}(\mu, d)}VWWI(T)$$ of vertex-weighted trees generated by the tuple $\langle\mu, d\rangle$, which minimize
the Wiener index.

\begin{definition}
The vertex-weighted tree $T$ \emph{induces} the tuple $\langle\mu, d\rangle$ on the vertex set $V = V(T)$ if $\mu(v) =
\mu_T(v), d(v) = d_T(v), v \in V(T)$. Clearly, the induced tuple $\langle\mu, d\rangle$ generates the tree $T$, i.e.,
$T \in \mathcal{T}(\mu,d)$.
\end{definition}

\begin{definition}
We will say that in the generating tuple $\langle\mu, d\rangle$
\emph{weights are degree-monotone} if for any pair of internal
vertices $m, m' \in M$ from $d(m) < d(m')$ it follows that $\mu(m)
\le \mu(m')$. We also require pendent vertices to have positive
weights: $d(v)=1\Rightarrow \mu(v)>0$.
\end{definition}

In this paper we show that if weights are degree-monotone in the tuple $\langle\mu, d\rangle$, then the set
$\mathcal{T}^*(\mu,d)$ consists of the trees built with the simple and efficient algorithm being a generalization of
the famous Huffman algorithm \cite{Huffman52} for construction of the binary tree of an optimal prefix code.

\subsection{Generalized Huffman Algorithm}

\begin{definition} A \emph{star} is a complete bipartite graph $K_{1,k}$, where a distinguished vertex, called a \emph{center}, is connected to
$k$ other vertices, called \emph{leaves}. For a star $S$, the set of its leaves is denoted with $L(S)$. It is clear
that $L(S)=W(S)$, except for the case of $S=K_{1,1}$, when $W(S)=V(S)$.
\end{definition}

\begin{definition}
Consider a generating tuple $\langle\mu, d\rangle$ with degree-monotone
weights. Let $m\in M$ be any internal vertex
having the least degree $d(m)$ among the vertices of the least weight in $M$,
i.e., $m\in \textrm{Argmin}\{d(u): u\in
\textrm{Argmin}_{v\in M}\mu(v)\}$. The \emph{minimal star} for the tuple
$\langle\mu, d\rangle$ is a vertex-weighted
star $S\in \mathcal{T}$ with the center $m$, $\mu_S(m)=\mu(m)$, and with
$d(m)-1$ leaves having $d(m)-1$ least weights
in $W$, i.e., $L(S)\subseteq W$, and $u\in L(S), v\in W\backslash L(S)
\Rightarrow \mu_S(u)=\mu(u)\le\mu(v)$. Denote
with $\underline{f}(\mu,d)$ the total weight of vertices of a minimal star.
\end{definition}

For a fixed tuple $\langle\mu,d\rangle$ the \emph{generalized Huffman algorithm} builds a tree $H \in
\mathcal{T}(\mu,d)$ as follows.

\textbf{Setup.} Define the vertex set $V_1 := V$ and the functions $\mu^1$ and $d^1$, which endow its vertices with
weights $\mu^1(v) := \mu(v)$ and degrees $d^1(v) := d(v)$, $v \in V_1$.

\textbf{Steps $i = 1, ..., q-1$.} Let the star $S_i$ be a minimal star for the tuple $\langle\mu^i,d^i\rangle$. Denote
its center with $m_i$. Define the set $V_{i+1} := V_i \backslash L(S_i)$ and functions $\mu^{i+1}, d^{i+1}$, endowing
its elements with weights and degrees as follows:
\begin{gather}
\mu^{i+1}(v) := \mu^i(v)\text{ for }v \neq m_i,  \mu^{i+1}(m_i) := \mu_{S_i} = \sum_{v\in V(S_i)}\mu^i(v),\nonumber\\
\label{eq_Huffman_tuples}d^{i+1}(v) := d^i(v)\text{ for }v \neq m_i, d^{i+1}(m_i) := 1.
\end{gather}

\textbf{Step $q$.} Consider a vertex $m_q \in M(\mu^q, d^q)$ (such a vertex is unique by construction), and a let $S_q$
be the star with the vertex set $V_q$ and the center $m_q$. We build a Huffman tree $H$ by setting $V(H) := V$, $E(H)
:= E(S_1) \cup ... \cup E(S_q)$, $\mu_H(v) := \mu(v)$, $v \in V$.

\begin{figure}[htpb]
\begin{center}
  \includegraphics[width=9cm]{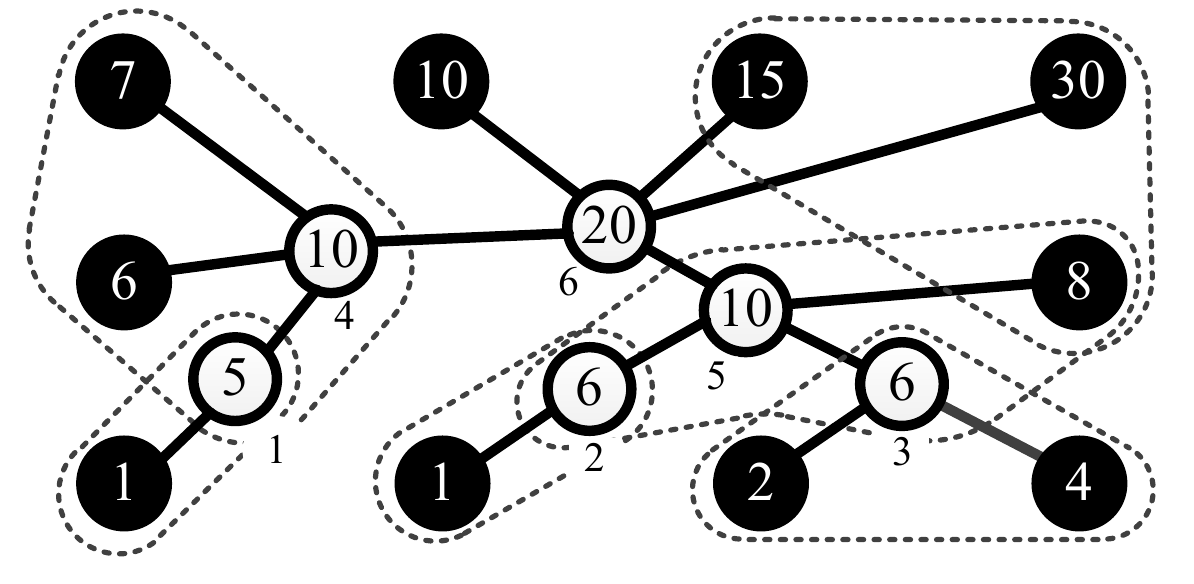}
  \caption{An example of Huffman tree construction}
  \label{fig_huffman}
\end{center}
\end{figure}

An example of Huffman tree construction is depicted in Fig.~\ref{fig_huffman}. Black circles correspond to pendent
vertices, numbers inside circles stand for vertex weights, those under circles show the order of star sequence centers.
All stars, except the last one, are surrounded by a dashed line.

Thus, the Huffman tree $H$ appears to be a union of minimal stars $S_1, ..., S_{q-1}$ for the corresponding generating
tuples and a ``finalizing'' star $S_q$. Below we refer to the sequence $S_1, ..., S_q$ as the \emph{star sequence} of a
Huffman tree $H$. In general, the Huffman tree is not unique, as more than one star sequence is possible. Let
$\mathcal{TH}(\mu,d)$ be the collection of Huffman trees generated by the tuple $\langle\mu, d\rangle$. The main result
of this paper can be stated as follows.

\begin{theorem}\label{theorem_main}
If weights are degree-monotone in a generating tuple $\langle\mu,d\rangle$, then $\mathcal{T}^*(\mu,d) =
\mathcal{TH}(\mu,d)$. In other words, only a Huffman tree minimizes the Wiener index over the set of trees whose
vertices have given weights and degrees.
\end{theorem}

In the following sections we prove auxiliary results, and return to the proof of Theorem~\ref{theorem_main} at the end
of Section~\ref{section_minimizing_index}.

\begin{figure}[htpb]
\begin{center}
  \includegraphics[width=12cm]{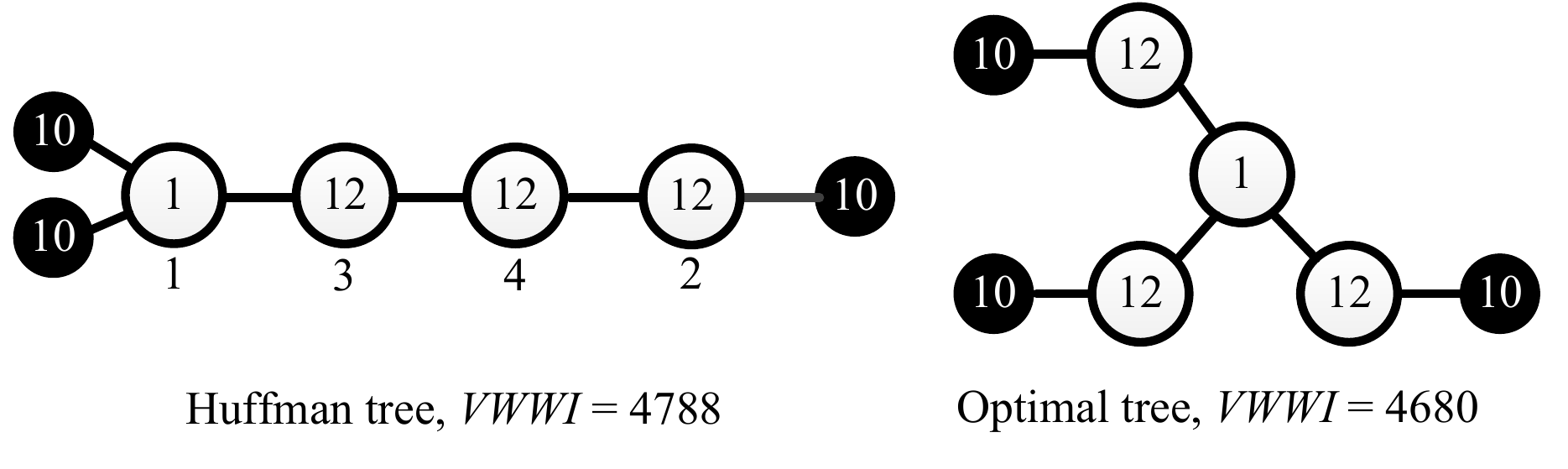}
  \caption{The counterexample for non-monotone weights}
  \label{fig_nonmonotonous}
\end{center}
\end{figure}

Please note that when $\mu(v)\equiv1$ for all $v\in V$, the Huffman tree becomes a ``greedy tree'' from \cite{Wang08}.
Fig.~\ref{fig_nonmonotonous} shows that weights' monotonicity is essential for Theorem \ref{theorem_main} (numbers
inside circles are vertex weights, those under circles show the order of star sequence centers).

\section{Properties of Huffman trees}\label{section_huffman_tree}
\subsection{Huffman algorithm for directed trees}

The index minimization problem becomes more tractable when studied for directed trees.

\begin{definition}
A (weighted) \emph{directed tree} is a connected directed graph with each vertex except the \emph{root }having the sole
outbound arc and the root having no outbound arcs.
\end{definition}

An arbitrary tree $T \in \mathcal{T}$ consisting of more than two vertices can be transformed into a directed tree
$T_r$ by choosing an \textbf{internal vertex} $r\in M(T)$ as a root, and replacing all its edges with arcs directed
towards the root. Let us denote with $\mathcal{R}$ the collection of all directed trees, which can be obtained in such
a way, and let $\mathcal{R}(\mu, d)$ stand for all directed trees obtained from $\mathcal{T}(\mu,d)$. Vice versa, in a
directed tree $T_r \in \mathcal{R}(\mu,d)$ replacing all arcs with edges makes some tree $T\in \mathcal{T}(\mu, d)$.

Let the arcs in a directed star be directed towards its center by definition.

If in a star sequence of a Huffman tree $H$ one replaces all stars with directed stars, then the union of the arcs of
these directed stars gives a \emph{directed Huffman tree }with the root at the center $m_q$ of the last star in the
sequence. Let $\mathcal{RH}(\mu,d)\subseteq \mathcal{R}(\mu,d)$ stand for the collection of directed Huffman trees
generated by the tuple $\langle\mu,d\rangle$.

\subsection{Vector of subordinate groups' weights and Wiener Index}

\begin{definition}
For an arbitrary vertex $v\in V(T)$ of the directed tree $T\in \mathcal{R}$ define its \emph{subordinate group}
$g_T(v)\subseteq V(T)$ as the set of vertices having the directed path to the vertex $v$ in the tree $T$ (the vertex
$v$ itself belongs to $g_T(v)$). The \emph{weight $f_T(v)$ of the subordinate group} $g_T(v)$ is defined as the total
vertex weight of the group: $f_T(v):=\sum_{u\in g_T(v)} \mu_T(u)$.
\end{definition}

In particular, all vertices in a directed tree $T\in \mathcal{R}$
are subordinated to its root $r$, i.e., $g_T(r)=V(T)$ and
$f_T(r)=\mu_T$. For example, if $T\in \mathcal{R}(\mu, d)$, then
$f_T(r)=\bar{\mu}$.

\begin{note}\label{note_positive}If all pendent vertices in $T$ have strictly positive weights,
then $f_T(v)>0$  for any $v\in V(T)$. In particular, it is true for
any $T\in \mathcal{R}(\mu, d)$, if weights in $\langle\mu,d\rangle$
are degree-monotone.
\end{note}

If some tree $T\in \mathcal{T}(\mu,d)$ is transformed into a directed tree $T_r\in \mathcal{R}(\mu,d)$ by choosing a
root $r$, the Wiener index can be written as \cite{Klavzar97, SkrekovskiGutman14}:
\begin{equation}\label{eq_VWWI_directed}
VWWI(T)=VWWI(T_r)=\sum_{v\in V\backslash \{r\}}f_{T_r}(v)(\bar{\mu}-f_{T_r}(v))=\sum_{v\in V\backslash
\{r\}}\chi(f_{T_r}(v)),
\end{equation}
where $\chi(x):=x(\bar{\mu}-x)$.

Equality (\ref{eq_VWWI_directed}) implies that all directed trees obtained from one tree $T\in \mathcal{T}(\mu,d)$
share the same value of the Wiener index. Thus, if we find the collection $\mathcal{R}^*(\mu,d):=\text{Argmin}_{T\in
\mathcal{R}(\mu, d)}VWWI(T)$ of directed trees minimizing the Wiener index, the collection $\mathcal{T}^*(\mu,d)$ is
obtained by replacing them with corresponding undirected trees.

As the root of a directed tree $T\in \mathcal{R}(\mu,d)$ is an
internal vertex, every pendent vertex has an outbound arc, so, for
every pendent vertex $w\in W$ in a directed tree $T\in
\mathcal{R}(\mu,d)$ $f_T(w)=\mu(w)$. Therefore, all directed trees
from $\mathcal{R}(\mu,d)$ enjoy the same weights of groups
subordinated to pendent vertices. Also, as noticed above,
$f_T(r)=\bar{\mu}$ for the root $r$ of any directed tree $T\in
\mathcal{R}(\mu,d)$. Thus, directed trees from $\mathcal{R}(\mu,d)$
differ only in the subordinate group weights of $q-1$ internal
vertices other than root.

\begin{definition} \cite{Marshall79,Zhang08}
For the real vector $\mathbf{x}=(x_1,..., x_p)$, $p\in \mathbb{N}$, denote with $\mathbf{x}_\uparrow=(x_{[1]},...,
x_{[p]})$ the vector, where all components of $\mathbf{x}$ are arranged in ascending order.
\end{definition}

\begin{definition}\label{def_vector}
For a directed tree $T\in \mathcal{R}(\mu,d)$ define a $(q-1)$-dimensional \emph{vector} $\mathbf{f}(T):=(f_T(m): m\in
M\backslash\{r\})_\uparrow$ \emph{of subordinate groups' weights}, where $r$ is the root of $T$.
\end{definition}

In the following proofs we combine the approach of \cite{Goubko04, Goubko06}, where Huffman tree has been proved to
minimize the \emph{sum of subordinate groups' weights} in case of zero-weighted internal vertices, and that by Zhang et
al \cite{Zhang08}, who minimized the Wiener index for unweighted trees having the given degree sequence.

\subsection{Basic Properties of Huffman Trees}

In Lemmas \ref{lemma_star_seq_weight}-\ref{lemma_Huffman_monotonicity} we consider a Huffman tree $H\in
\mathcal{RH}(\mu,d)$ with a star sequence  $S_1, ..., S_q$, and vertices $m_1, ..., m_q$ being the centers of stars
$S_1, ..., S_q$ respectively.

\begin{lemma}\label{lemma_star_seq_weight} $\mathbf{f}(H) = (\mu^2(m_1), \mu^3(m_2),..., \mu^q(m_{q-1})) = (\underline{f}(\mu^1, d^1), ...,
\underline{f}(\mu^{q-1}, d^{q-1}))$, where tuples $\langle\mu^i, d^i\rangle, i = 1, ..., q - 1$, are defined by formula
\emph{(\ref{eq_Huffman_tuples})}.
\begin{proof}
The definition of a minimal star implies that $\mu_{S_i}=\underline{f}(\mu^i, d^i)$. By construction of tuples
$\langle\mu^i,d^i\rangle$ we have $f_H(m_i)=\sum_{v\in V(S_i)}\mu_i(v)=\mu^{i+1}(m_i)$, and thus,
$f_H(m_i)=\mu^{i+1}(m_i)=\underline{f}(\mu^i, d^i), i=1,...,q-1$. One can easily see that $\underline{f}(\mu^i,d^i)\le
\underline{f}(\mu^{i+1},d^{i+1}), i=1,...,q-2$, from which the statement of the lemma follows immediately.
\end{proof}
\end{lemma}

\begin{lemma}\label{lemma_stars_monotonicity}
From $v\in L(S_i), v'\in L(S_j)$, and $i<j$ it follows that $f_H(v)\le f_H(v')$.
\begin{proof}
Suppose, by contradiction, that $f_H(v)> f_H(v')$. As $v\in L(S_i), v'\in L(S_j)$, and $i<j$, a vertex $v''\in g_H(v')$
exists, which also belongs to $W(\mu^i,d^i)$ (otherwise the vertex $v'$ cannot belong to the set $W(\mu^j, d^j)$, as
the tuple $\langle\mu^j, d^j\rangle$ is defined later, at the $(j-1)$-th step of the algorithm).

By definition of a subordinate group, a path exists from the vertex
$v''$ to $v'$ in $H$, which immediately implies that $f_H(v') \ge
f_H(v'')$, and, by assumption, $f_H(v) > f_H(v'')$. Then the vertex
$v$ cannot be a leaf of $S_i$ by definition of a minimal star, as
the set $W(\mu^i,d^i)$ contains the vertex $v''$, which does not
belong to the minimal star $S_i$, but has the weight $\mu^i(v'') <
\mu^i(v)$ (since, by Lemma \ref{lemma_star_seq_weight}, $\mu^i(v'')
= f_H(v''), \mu^i(v) = f_H(v))$. We obtain a contradiction, so the
lemma is correct.
\end{proof}
\end{lemma}

\begin{lemma}\label{lemma_Huffman_monotonicity}
If weights are degree-monotone in $\langle\mu,d\rangle$, then for
any $H\in \mathcal{RH}(\mu,d)$
\begin{equation}\label{eq_weights_monotone}
[vm,v'm'\in E(H), m \neq m', f_H(v)<f_H(v')] \Rightarrow f_H(m) < f_H(m').
\end{equation}
\begin{proof}
Suppose, by contradiction, that a pair of arcs $vm_i, v'm_j\in E(H)$ exists, such that $m_i \neq m_j$, $f_H(v) <
f_H(v')$, but $f_H(m_i) \ge f_H(m_j)$. In case of strict inequality $f_H(m_i) > f_H(m_j)$, from Lemma
\ref{lemma_star_seq_weight}, it follows that $i > j$. Then (since presence of the arcs $vm_i$ and $v'm_j$ implies that
$v \in L(S_i)$ and $v' \in L(S_j)$), by Lemma \ref{lemma_stars_monotonicity}, $f_H(v) \ge f_H(v')$. We obtain a
contradiction, and, since $m_i \neq m_j$, we are left with the sole case of $i < j$ and $f_H(m_i) = f_H(m_j)$.

Since, by Lemma \ref{lemma_stars_monotonicity}, for every pair of
vertices $u \in L(S_i)$, $u' \in L(S_j)$ we have $f_H(u) \le
f_H(u')$, and, by construction of the Huffman tree, $\mu(m_i) \le
\mu(m_j)$, and also, from degree-monotonicity of weights in
$\langle\mu,d\rangle$ we have $d(m_i)\le d(m_j)$, the equality
$f_H(m_i) = f_H(m_j)$ is possible only if $\mu(m_i) = \mu(m_j)$, and
$f_H(u) = f_H(u')$ for all $u \in L(S_i), u' \in L(S_j)$. However,
by assumption, $v \in L(S_i), v' \in L(S_j)$ and $f_H(v) < f_H(v')$.
The obtained contradiction completes the proof.
\end{proof}
\end{lemma}

\subsection{Vector of subordinate groups' weights in Huffman Trees}
In this paragraph we show that all directed Huffman trees share the same vector of subordinate groups' weights, and no
other tree enjoys this vector of subordinate groups' weights. These results allow us to move the index minimization
problem into the space of vectors of subordinate groups' weights for    directed trees from $\mathcal{R}(\mu,d)$.

\begin{definition}
Consider a tuple $\langle\mu, d\rangle$ of functions (not necessarily the generating one) defined on the set $V$, and a
tuple $\langle\mu', d'\rangle$ defined on the set $V'$. A bijection $\sigma: V \rightarrow V'$ \emph{preserves weights
and degrees} if $\mu(v) = \mu'(\sigma(v))$, $d(v) = d'(\sigma(v))$, $v \in V$.
\end{definition}

\begin{lemma}\label{lemma_Huffman_bijection1}
Consider a generating tuple $\langle\mu, d\rangle$ on the set $V$, a tuple $\langle\mu', d'\rangle$ on the set $V'$,
and a bijection $\sigma: V \rightarrow V'$ preserving weights and degrees. If $H \in \mathcal{RH}(\mu,d)$ is a directed
Huffman tree, then there exists a Huffman tree $H' \in \mathcal{RH}(\mu',d')$ such that $\mathbf{f}(H) =
\mathbf{f}(H')$.
\begin{proof}
Consider a star sequence $S_1,...,S_q$ of the Huffman tree $H$, with
$m_1, ..., m_q$ being the centers of stars $S_1,...,S_q$
respectively. The Huffman algorithm takes care only of vertex
weights and degrees, so, replacing all vertices in stars
$S_1,...,S_q$ with their images under the bijection $\sigma(\cdot)$,
we obtain the sequence $\sigma(S_1),...,\sigma(S_q)$ of stars, which
give some Huffman tree $H'\in \mathcal{RH}(\mu',d')$ as their union.

As the group $g_{H'}(\sigma(m_i))$ subordinated in the directed tree
$H'$ to the image $\sigma(m_i)$ of the vertex $m_i$ coincides with
the image $\sigma(g_H(m_i))$ of the subordinate group of the vertex
$m_i$ in the directed tree $H$, we obtain $f_H(m_i) =
f_{H'}(\sigma(m_i))$. So, according to Definition \ref{def_vector},
$\mathbf{f}(H) = \mathbf{f}(H')$.
\end{proof}
\end{lemma}

\begin{lemma}\label{lemma_Huffman_bijection2}
If $S$ and $S'$ are two different minimal stars for the tuple
$\langle\mu, d\rangle$, then a bijection $\sigma: L(S)\rightarrow
L(S')$ preserving weights and degrees can be established between the
leaf sets $L(S)$ and $L(S')$ of these stars.
\begin{proof}
By definition of a minimal star, sets $L(S)$ and $L(S')$ consist of the same
number of elements. Define the vectors $\mathbf{w} := (\mu(v): v \in
L(S))_\uparrow$ and $\mathbf{w'} := (\mu(v): v \in
L(S'))_\uparrow$. Since both $L(S)$ and $L(S')$ include the same number of vertices having the minimum weight in $W$,
it is clear that $\mathbf{w} = \mathbf{w'}$. The desired bijection is built by matching sequentially vertices inducing
the first, the second, etc, components of the vectors $\mathbf{w}$ and $\mathbf{w'}$.
\end{proof}
\end{lemma}

\begin{definition}
A directed star $S$ with the center $m \in M(T)$ is called the \emph{lower star} of a directed tree $T \in
\mathcal{R}$, if $V(S) = g_T(m)$ and $\mu_S(v) = \mu_T(v), v \in V(S)$.
\end{definition}

\begin{definition}\label{def_rollup}
Let $m \in M(T)$ be an internal vertex in a directed tree $T \in \mathcal{R}$. The \emph{$m$-rollup} of $T$ is a
directed tree $\underline{T} \in \mathcal{R}$ obtained from $T$ by deleting the set of vertices $g_T(m)\backslash\{m\}$
along with their incident arcs, and setting $\mu_{\underline{T}}(m):=f_T(m)$.
\end{definition}

Please note that if a directed tree $R$ is a contraction of $T$ to the vertex set $V(R):=g_T(m)$, and $m$ is not a root
of $T$, then $\mathbf{f}(T) = (\mathbf{f}(R), f_T(m), \mathbf{f}(\underline{T}))_\uparrow$.

\begin{lemma}\label{lemma_Huffman_rollup}
Consider the star sequence $S_1,...,S_q$ of a Huffman tree $H\in \mathcal{RH}(\mu,d)$ with the vertex $m_1$ being the
center of the star $S_1$. If the tuple $\langle\mu',d'\rangle$ is induced by the $m_1$-rollup $\underline{H}$ of the
Huffman tree $H$, then $\underline{H}\in \mathcal{RH}(\mu', d')$. In other words, the $m_1$-rollup of a Huffman tree
appears to be a Huffman tree for the induced generating tuple.
\begin{proof}
By construction of the Huffman tree the tuple $\langle\mu', d'\rangle$ coincides with the tuple $\langle\mu^2,
d^2\rangle$ from the Huffman algorithm. Thus, $S_2$ is a minimal star for $\langle\mu', d'\rangle$, which implies that
the stars $S_3, ..., S_q$ are minimal stars for the corresponding generating tuples defined with formula
(\ref{eq_Huffman_tuples}). As $E(\underline{H}) = E(S_2) \cup...\cup E(S_q)$, by definition of a Huffman tree we obtain
$\underline{H} \in \mathcal{RH}(\mu', d')$.
\end{proof}
\end{lemma}

\begin{lemma}\label{lemma_all_Huffman_same_weights}
All Huffman trees share the same vector of subordinate groups' weights, i.e., if $T, H \in \mathcal{RH}(\mu, d)$, then
$\mathbf{f}(T) = \mathbf{f}(H)$.
\begin{proof}
Employ induction on the number of internal vertices $q$. For $q=1$ the vector of subordinate groups' weights has zero
components, thus, the lemma obviously holds.

Suppose the lemma holds for all $q' < q$. Let us prove that it also holds for the set $V$ with $q$ internal vertices.
Denote $\mathbf{f}(T) = (f_1, ..., f_{q-1}), \mathbf{f}(H) = (f_1', ..., f_{q-1}')$. On the first step of the Huffman
algorithm some minimal stars $S$ and $S'$ with the centers $m$ and $m'$ are added to the trees $T$ and $H$
respectively, thus, $f_1 = f_1' = \underline{f}(\mu, d)$. Consider the $m_1$-rollup $\underline{T}$ of the tree $T$ and
the $m_1'$-rollup $\underline{H}$ of the tree $H$. Let $\underline{T}$ induce the tuple $\langle\mu', d'\rangle$ and
$\underline{H}$ induce the tuple $\langle\mu'', d''\rangle$. From Lemma \ref{lemma_Huffman_rollup}, $\underline{T} \in
\mathcal{RH}(\mu', d')$, $\underline{H} \in \mathcal{RH}(\mu'', d'')$. By Lemma \ref{lemma_star_seq_weight},
$\mathbf{f}(\underline{T}) = (f_2, ..., f_{q-1})$, $\mathbf{f}(\underline{H}) = (f_2', ..., f_{q-1}')$.

From Lemma \ref{lemma_Huffman_bijection2}, a bijection can be established between elements of the sets $L(S)$ and
$L(S')$, which preserves weights and degrees. So, obviously, an analogous bijection $\sigma$ can be established between
the elements of the residual sets $V(\underline{T}) = V\backslash L(S)$ (with the generating tuple $\langle\mu', d
'\rangle$) and $V(\underline{H}) = V\backslash L(S')$ (with the generating tuple $\langle\mu'', d''\rangle$), which
also preserves weights and degrees. Thus, by Lemma \ref{lemma_Huffman_bijection1}, there exists such a Huffman tree
$\underline{\underline{H}} \in \mathcal{RH}(\mu', d')$ that $\mathbf{f}(\underline{\underline{H}}) =
\mathbf{f}(\underline{H})$.

There are $q-1$ internal vertices in the tree $\underline{T}$, so, by inductive assumption $(f_2, ..., f_{q-1}) =
(f_2', ..., f_{q-1}')$ and, since $f_1 = f_1'$, the proof is complete.
\end{proof}
\end{lemma}

\begin{lemma}\label{lemma_same_weights_is_Huffman}
If a tree has the same vector of subordinate groups' weights as some Huffman tree, it has to be a Huffman tree itself.
In other words, for $H \in \mathcal{RH}(\mu,d), T \in \mathcal{R}(\mu,d)$ from $\mathbf{f}(H) = \mathbf{f}(T)$ it
follows that $T \in \mathcal{RH}(\mu,d)$.
\begin{proof}
We again employ induction on the number of internal vertices $q$. For $q=1$ the vector of subordinate groups' weights
has zero components, but $H = T$, since the collection $\mathcal{R}(\mu,d)$ consists of the sole directed tree (the
star). Assume the lemma is valid for all $q' < q$; let us prove that it also holds for the vertex set $V$ with $q$
internal vertices.

Denote for short $\mathbf{f}(H) = \mathbf{f}(T) = (f_1, ..., f_{q-1})$. By construction of the Huffman tree $H$, $f_1 =
\underline{f}(\mu, d)$. Every star with the total vertex weight $\underline{f}(\mu, d)$ is minimal, so, some minimal
star $S_1$ for the tuple $\langle\mu, d\rangle$ must be a part of the tree $T$; $H$ contains some minimal star $S_1'$
by definition. Denote with $m_1, m_1'$ respectively the centers of these stars.

Let the tuple $\langle\mu', d'\rangle$ be induced by the $m_1$-rollup $\underline{T}$ of the directed tree $T$, and the
tuple $\langle\mu'', d''\rangle$ be induced by the $m_1'$-rollup $\underline{H}$ of the directed Huffman tree $H$. By
Lemma \ref{lemma_Huffman_rollup}, $\underline{H} \in \mathcal{RH}(\mu'', d'')$. Moreover, by Lemma
\ref{lemma_star_seq_weight}, $\mathbf{f}(\underline{T}) = \mathbf{f}(\underline{H}) = (f_2, ..., f_{q-1})$.

By analogy with the proof of Lemma \ref{lemma_all_Huffman_same_weights}, between the vertex sets $V(\underline{T})$
(with the tuple $\langle\mu', d'\rangle$) and $V(\underline{H})$ (with the tuple $\langle\mu'', d''\rangle$) one can
establish a bijection $\sigma$ preserving weights and degrees, so, by Lemma \ref{lemma_Huffman_bijection1}, such a
Huffman tree $\underline{\underline{H}} \in \mathcal{RH}(\mu', d')$ exists that $\mathbf{f}(\underline{H}) =
\mathbf{f}(\underline{\underline{H}})$. Then we have $\mathbf{f}(\underline{\underline{H}}) = \mathbf{f}(\underline{T})
= (f_2, ..., f_{q-1})$, and, by inductive assumption, $\underline{T}$ is a Huffman tree for the tuple $\langle\mu',
d'\rangle$. Let $S_2, ..., S_q$ be its star sequence. Then the tree $T$ can be obtained as a union of $\underline{T}$
and the minimal star $S_1$, and, thus, $T \in \mathcal{RH}(\mu, d)$.
\end{proof}
\end{lemma}

To sum up, Lemmas \ref{lemma_all_Huffman_same_weights} and \ref{lemma_same_weights_is_Huffman} say that if some Huffman
tree $H$ has the vector $\mathbf{f}(H)$ of subordinate groups' weights, then all Huffman trees, and only they, have
this vector of subordinate groups' weights.

\begin{corollary}
If $H, H' \in \mathcal{RH}(\mu, d)$ are two directed Huffman trees, then $VWWI(H) = VWWI(H')$.
\begin{proof}
From equation (\ref{eq_VWWI_directed}) we know that the value of the index is determined by the components of vectors
$\mathbf{f}(H), \mathbf{f}(H')$, and also by the weights of pendent vertices of trees $H$ and $H'$. From Lemma
\ref{lemma_all_Huffman_same_weights} we learn that $\mathbf{f}(H) = \mathbf{f}(H')$, so, since the trees $H$ and $H'$
enjoy the same weights of pendent vertices, we induce that the index has the same value for both trees.
\end{proof}
\end{corollary}

Therefore, to justify Theorem \ref{theorem_main} it is enough to prove that the vector of subordinate groups' weights
originated from some Huffman tree minimizes \emph{VWWI} over all directed trees in the collection $\mathcal{R}(\mu,d)$.
We postpone the proofs to the next section.

\section{Huffman Trees and Majorization}\label{section_majorization}

\subsection{Notion of Vectors' Majorization}

Let us recall that notation $\mathbf{x}_\uparrow = (x_{[1]}, ..., x_{[p]})$ stands for the vector where all components
of a real vector $\mathbf{x} = (x_1, ..., x_p), p \in \mathbb{N}$, are arranged in the ascending order.

\begin{definition} \textbf{\cite{Marshall79,Zhang08}}
A non-negative vector $\mathbf{x} = (x_1, ..., x_p)$, $p \in
\mathbb{N}$, \emph{weakly majorizes} a non-negative vector
$\mathbf{y} = (y_1, ..., y_p)$ (which is denoted with
$\mathbf{y}\preceq^w \mathbf{x}$ or $\mathbf{x}\succeq^w
\mathbf{y}$) if
$$\sum_{i=1}^k x_{[i]} \le \sum_{i=1}^k y_{[i]} \text{ for all }k=1,...,p.$$
Moreover, if $\mathbf{x}_\uparrow \neq \mathbf{y}_\uparrow$, then
$\mathbf{x}$ is said to \emph{strictly weakly majorize} $\mathbf{y}$
(which is denoted with $\mathbf{y}\prec^w \mathbf{x}$ or
$\mathbf{x}\succ^w \mathbf{y}$).
\end{definition}

We will need the following properties of weak majorization.

\begin{lemma}\label{lemma_Zhang_b}\textbf{\emph{\cite{Marshall79,Zhang08}}}
Consider a positive number $b>0$ and two non-negative vectors,
$\mathbf{x} = (x_1, ..., x_k, y_1, ..., y_l)$ and $\mathbf{y} = (x_1
+ b, ..., x_k + b, y_1 - b, ..., y_l - b)$, such that $0 \le k \le
l$. If $x_i \ge y_i$ for $i = 1, ..., k$, then $\mathbf{x}\prec^w
\mathbf{y}$.
\end{lemma}

\begin{lemma}\label{lemma_Zhang_xy}\textbf{\emph{\cite{Marshall79,Zhang08}}}
If $\mathbf{x}\preceq^w \mathbf{y}$ and $\mathbf{x}' \preceq^w
\mathbf{y'}$, then $(\mathbf{x},\mathbf{x'})\preceq^w
(\mathbf{y},\mathbf{y'})$, where $(\mathbf{x},\mathbf{x'})$ means
concatenation of vectors $\mathbf{x}$ and $\mathbf{x'}$. Moreover,
if $\mathbf{x}' \prec^w \mathbf{y'}$, then
$(\mathbf{x},\mathbf{x'})\prec^w (\mathbf{y},\mathbf{y'})$.
\end{lemma}

\begin{lemma}\label{lemma_Zhang_concave}\textbf{\emph{\cite{Marshall79,Zhang08}}}
If $\chi(x)$ is an increasing concave function, and
$(x_1,...,x_p)\preceq^w (y_1,...,y_p)$, then
$\sum_{i=1}^p\chi(x_i)\ge \sum_{i=1}^p\chi(y_i)$, and equality is
possible only when $(x_1,...,x_p)_\uparrow=(y_1,...,y_p)_\uparrow$.
\end{lemma}

\subsection{Transformations of Trees and Majorization}

The following lemmas play the same role in our proofs as Lemmas 3.1-3.5 in \cite{Zhang08}. Some novelty is originated
from accounting for variations in internal vertex weights.

\begin{lemma}\label{lemma_major1}\sloppy
Suppose a directed tree $T\in \mathcal{R}(\mu, d)$ contains the
disjoint paths $(v, m_1, ..., m_k, m)$ and $(v', m_1', ..., m_l',
m)$ from vertices $v, v' \in V$ to some vertex $m \in M$, and
suppose that $1 \le k \le l$, $f_T(v) < f_T(v')$, $f_T(m_i) \ge
f_T(m_i')$, $i = 1, ..., k$. If the directed tree $T'$ is obtained
from $T$ by deleting the arcs $vm_1, v'm_1'$ and adding the arcs
$v'm_1$ and $vm_1'$ instead, then $T'\in \mathcal{R}(\mu, d)$ and
$\mathbf{f}(T')\succ^w \mathbf{f}(T)$.
\begin{proof}
Clearly, $T'\in\mathcal{R}(\mu, d)$, since vertex degrees and weights do not change during the transformation. Denote
$b := f_T(v') - f_T(v)> 0$. In the tree $T'$ weights of the groups subordinated to the vertices $m_1, ..., m_k \in M$
increase by $b$ (i.e., $f_{T'}(m_i) = f_T(m_i) + b, i = 1, ..., k$), weights of the groups subordinated to the vertices
$m_1', ..., m_l' \in M$ decrease by $b$ (i.e., $f_{T'}(m_i') = f_T(m_i') - b, i = 1, ..., l$), weights of all other
vertices (including $m$) do not change. Therefore, by Lemma \ref{lemma_Zhang_b},
$$\mathbf{y} := (f_{T'}(m_1), ..., f_{T'}(m_k), f_{T'}(m_1'), ..., f_{T'}(m_l'))= $$
$$=(f_T(m_1) + b, ..., f_T (m_k) + b, f_T(m_1') - b, ..., f_T (m_l') - b) \succ^w $$
$$\succ^w (f_T (m_1), ..., f_T (m_k), f_T (m_1'), ..., f_T (m_l')) =: \mathbf{x}.$$

If one denotes with $\mathbf{z}$ the vector of (unchanged) weights
of groups subordinated to all other internal vertices of $T$
distinct from the root, then, by Lemma \ref{lemma_Zhang_xy},
$\mathbf{f}(T')=(\mathbf{y},\mathbf{z})\succ^w
(\mathbf{x},\mathbf{z})=\mathbf{f}(T)$.
\end{proof}\end{lemma}

\begin{lemma}\label{lemma_major2}\sloppy
Consider a directed tree $T\in \mathcal{R}(\mu, d)$ containing the
disjoint paths $(v, m_1, ..., m_k, m)$ and $(v', m_1', ..., m_l',
m)$ from vertices $v, v' \in V$ to some vertex $m \in M$, and
suppose that $1 \le l \le k$, $f_T(v) < f_T(v')$, $f_T(m_1) =
f_T(m_1')$, and $f_T(m_i) \le f_T(m_i')$, $i = 2, ..., l$. Then such
a directed tree $T'\in \mathcal{R}(\mu, d)$ exists that
$\mathbf{f}(T')\succ^w \mathbf{f}(T)$.
\begin{proof}
Introduce the notation
$$u=\begin{cases} m_2, &\mbox{if } k\ge2 \\
m, & \mbox{if } k=1,\end{cases}\hspace{20pt}
u'=\begin{cases} m_2', &\mbox{if } l\ge2 \\
m, & \mbox{if } l=1,\end{cases}$$ and consider the tree $T'$ obtained from $T$ by deleting the arcs $vm_1, v'm_1',
m_1u, m_1'u'$ and adding the arcs $v'm_1, vm_1', m_1u'$, and $m_1'u$ instead. We have $T'\in\mathcal{R}(\mu, d)$, since
vertex degrees and weights do not change during the transformation.

Denote $b := f_T(v') - f_T(v)> 0$. In the tree $T'$ weights of the groups subordinated to the vertices $m_1', m_2, ...,
m_k \in M$ decrease by $b$ (i.e., $f_{T'}(m_1') = f_T(m_1') - b$, $f_{T'}(m_i) = f_T(m_i) - b, i = 2, ..., k$), weights
of the groups subordinated to the vertices $m_1, m_2', ..., m_l' \in M$ increase by $b$ (i.e., $f_{T'}(m_1) = f_T(m_1)
+ b$, $f_{T'}(m_i') = f_T(m_i') + b, i = 2, ..., l$), while weights of all other vertices (including $m$) do not
change. Therefore, by Lemma \ref{lemma_Zhang_b},
$$\mathbf{y} := (f_{T'}(m_1'), f_{T'}(m_2), ..., f_{T'}(m_k), f_{T'}(m_1), f_{T'}(m_2'), ..., f_{T'}(m_l'))=$$
$$=(f_T(m_1') - b, f_T(m_2) - b, ..., f_T (m_k) - b, f_T(m_1) + b, f_T(m_2') + b, ..., f_T (m_l') + b) \succ^w$$
$$\succ^w (f_T (m_1), ..., f_T (m_k), f_T (m_1'), ..., f_T (m_l')) =: \mathbf{x}.$$

If one denotes with $\mathbf{z}$ the vector of weights of groups
subordinated to all other internal vertices of $T$ distinct from the
root, then, by Lemma \ref{lemma_Zhang_xy},
$\mathbf{f}(T')=(\mathbf{y},\mathbf{z})\succ^w
(\mathbf{x},\mathbf{z})=\mathbf{f}(T)$.
\end{proof}
\end{lemma}

\begin{lemma}\label{lemma_major3}\sloppy
Consider a directed tree $T\in \mathcal{R}(\mu, d)$, which contains
the paths $(v, m)$ and $(v', m_1', ..., m_l', m)$ from vertices $v,
v' \in V$ to some vertex $m \in M$, and suppose that $l \ge 1$ and
$f_T(v) < f_T(v')$. If the directed tree $T'$ is obtained from $T$
by deleting the arcs $vm, v'm_1'$ and adding the arcs $v'm$ and
$vm_1'$ instead, then $T'\in \mathcal{R}(\mu, d)$ and
$\mathbf{f}(T')\succ^w \mathbf{f}(T)$.
\begin{proof}
Since vertex degrees and weights do not change during the transformation, $T'\in\mathcal{R}(\mu, d)$. Denote $b :=
f_T(v') - f_T(v)> 0$. In the tree $T'$ weights of the groups subordinated to the vertices $m_1', ..., m_l' \in M$
decrease by $b$ (i.e., $f_{T'}(m_i') = f_T(m_i') - b, i = 1, ..., l$), while weights of all other vertices do not
change. Therefore, by Lemma \ref{lemma_Zhang_b},
$$\mathbf{y} := (f_{T'}(m_1'), ..., f_{T'}(m_l'))=(f_T(m_1') - b, ..., f_T (m_l') - b) \succ^w (f_T (m_1'), ..., f_T
(m_l')) =: \mathbf{x}.$$

If $\mathbf{z}$ is the vector of weights of groups subordinated to
all other internal vertices of $T$ distinct from the root, then, by
Lemma \ref{lemma_Zhang_xy},
$\mathbf{f}(T')=(\mathbf{y},\mathbf{z})\succ^w
(\mathbf{x},\mathbf{z})=\mathbf{f}(T)$.
\end{proof}
\end{lemma}

\begin{lemma}\label{lemma_major4}\sloppy
Suppose weights are degree-monotone in a generating tuple
$\langle\mu,d\rangle$ and consider a directed tree $T\in
\mathcal{R}(\mu, d)$ containing the disjoint paths $(v, m_1, ...,
m_k, m)$ and $(v', m_1', ..., m_l', m)$ from vertices $v, v' \in M$
to some vertex $m \in M$. Suppose that $0 \le k \le l$, $d_T(v') -
d_T(v) = \Delta > 0$, $f_T(v) \ge f_T(v')$, $f_T(m_i) \ge
f_T(m_i')$, $i = 1, ..., k$. Then there exists a directed tree
$T'\in \mathcal{R}(\mu, d)$ such that $\mathbf{f}(T')\succ^w
\mathbf{f}(T)$.
\begin{proof}
Let the vertex $v$ have $d > 0$ inbound arcs in $T$. Introduce the notation
$$u=\begin{cases} m_1, &\mbox{if } k\ge1 \\
m, & \mbox{if } k=0,\end{cases}\hspace{20pt}
u'=\begin{cases} m_1', &\mbox{if } l\ge1 \\
m, & \mbox{if } l=0,\end{cases}$$ and consider the tree $T'$ obtained from $T$ by replacing the arcs $vu, v'u'$  with
the arcs $v'u, vu'$, redirecting all $d$ inbound arcs of the vertex $v$ in $T$ to the vertex $v'$, and redirecting
arbitrary $d$ inbound arcs of the vertex $v'$ in $T$ to the vertex $v$. We have $T'\in\mathcal{R}(\mu, d)$, since
vertex degrees and weights do not change during the transformation.

Let $u_1, ..., u_\Delta$ be those $\Delta$ vertices, for which
outbound arcs to the vertex $v'$ in the tree $T$ survived in the
tree $T'$, and introduce $b := f_T(u_1) + ... + f_T(u_\Delta) +
[\mu(v') - \mu(v)]$. Weights are degree-monotone in
$\langle\mu,d\rangle$, so we have $\mu(v') - \mu(v) \ge 0$. Since
$\Delta> 0$, from Note \ref{note_positive} it follows that $b
> 0$.

In the tree $T'$ weights of the groups subordinated to the vertices $m_1, ..., m_k$ (when $k>0$) increase by $b$,
weights of the groups subordinated to the vertices $m_1', ..., m_l'$ (when $l>0$) decrease by $b$. Also we have
$f_{T'}(v') - f_T(v) = -[f_{T'}(v) - f_T(v')] = b$. Weights of all other vertices (including $m$) do not change.
Therefore, by Lemma \ref{lemma_Zhang_b},
$$\mathbf{y} := (f_{T'}(v'), f_{T'}(m_1), ..., f_{T'}(m_k), f_{T'}(v), f_{T'}(m_1'), ..., f_{T'}(m_l')) = $$
$$= (f_T(v) + b, f_T(m_1) + b, ..., f_T(m_k) + b, f_T(v') - b, f_T(m_1') - b, ..., f_T(m_l') - b) \succ^w $$
$$\succ^w  (f_T(v), f_T(m_1), ..., f_T(m_k), f_T(v'), f_T(m_1'), ..., f_T(m_l')) =: \mathbf{x}.$$

If $\mathbf{z}$ is the vector of weights of groups subordinated to
all other internal vertices of $T$ distinct from the root, then, by
Lemma \ref{lemma_Zhang_xy},
$\mathbf{f}(T')=(\mathbf{y},\mathbf{z})\succ^w
(\mathbf{x},\mathbf{z})=\mathbf{f}(T)$.
\end{proof}
\end{lemma}

\begin{lemma}\label{lemma_major5}\sloppy
Consider a directed tree $T\in \mathcal{R}(\mu, d)$ containing the disjoint paths $(v, m_1, ..., m_k, m)$ and $(v',
m_1', ..., m_l', m)$ from vertices $v, v' \in M$ to some vertex $m \in M$, and suppose that $0\le k \le
 l$, $d(v) = d(v')$, $\mu(v) < \mu(v')$, $f_T(v) \ge f_T(v')$, and $f_T(m_i) \ge f_T(m_i')$, $i = 1, ..., k$.
 If the directed tree $T'$ is obtained from $T$ by swapping all incident arcs of vertices $v$ and $v'$, then
 $T'\in \mathcal{R}(\mu, d)$ and $\mathbf{f}(T')\succ^w \mathbf{f}(T)$.
\begin{proof}
It is clear that $T'\in\mathcal{R}(\mu, d)$. Denote $b := \mu(v') - \mu(v) > 0$. The rest of the proof repeats the one
of Lemma \ref{lemma_major4}.
\end{proof}
\end{lemma}

\begin{lemma}\label{lemma_major5a}\sloppy
Consider a path $(v', m_1, ..., m_k, v)$, $k\ge0$, in a directed
tree $T\in \mathcal{R}(\mu, d)$, and suppose that $d(v) = d(v')$,
$\mu(v') > \mu(v)$. If the directed tree $T'$ is obtained from $T$
by swapping all incident arcs of vertices $v$ and $v'$, then $T'\in
\mathcal{R}(\mu, d)$ and $\mathbf{f}(T')\succ^w \mathbf{f}(T)$.
\begin{proof}
It is clear that $T'\in\mathcal{R}(\mu, d)$. Denote $b := \mu(v') - \mu(v) > 0$. Then $f_{T'}(v')=f_T(v)$,
$f_{T'}(v)=f_T(v')-b$, $f_{T'}(m_i)=f_T(m_i)-b$, $i=1,...,k$. Weights of all other vertices do not change. Therefore,
by Lemmas \ref{lemma_Zhang_b} and \ref{lemma_Zhang_xy},
$$\mathbf{y} := (f_{T'}(v), f_{T'}(m_1), ..., f_{T'}(m_k), f_{T'}(v')) = $$
$$= (f_T(v') - b, f_T(m_1) - b, ..., f_T(m_k) - b, f_T(v)) \succ^w $$
$$\succ^w  (f_T(v'), f_T(m_1), ..., f_T(m_k), f_T(v)) =: \mathbf{x}.$$

If $\mathbf{z}$ is the vector of weights of groups subordinated to
all other internal vertices of $T$ distinct from the root, then, by
Lemma \ref{lemma_Zhang_xy},
$\mathbf{f}(T')=(\mathbf{y},\mathbf{z})\succ^w
(\mathbf{x},\mathbf{z})=\mathbf{f}(T)$.
\end{proof}
\end{lemma}

\begin{lemma}\label{lemma_major6}\sloppy
Suppose weights are degree-monotone in a generating tuple $\langle\mu,d\rangle$ and consider a directed tree $T\in
\mathcal{R}(\mu, d)$. Let $T$ contain an arc $mm'\in E(T)$, and suppose that $d_T(m) - d_T(m') = \Delta
> 0$. Then such a directed tree $T'\in
\mathcal{R}(\mu, d)$ exists that $\mathbf{f}(T')\succ^w
\mathbf{f}(T)$.
\begin{proof}
If the vertex $m'$ has an outbound arc in the tree $T$, denote this arc with $m'u$. Let the vertex $m'$ have $d \ge 0$
inbound arcs from the vertices other than $m$. Consider a directed tree $T'$ obtained from $T$ by replacing the arc
$mm'$ with the inverse arc $m'm$, replacing the arc $m'u$ (if it presents) with the arc $mu$, redirecting to the vertex
$m$ all $d$ arcs entering the vertex $m'$ from the vertices other than $m$ in $T$, and redirecting to the vertex $m'$
as many (arbitrary) inbound arcs of the vertex $m$ in $T$ as needed to restore its degree $d(m')$ (we are enough to
redirect $d$ arcs in case of $m'$ being a root in $T$, and $d+1$ arcs otherwise). Since vertex degrees and weights do
not change during the transformation, we have $T'\in\mathcal{R}(\mu, d)$.

Let $u_1, ..., u_\Delta$ be those $\Delta$ vertices, for which
outbound arcs to the vertex $m$ in the tree $T$ survived in the tree
$T'$, and introduce $b := f_T(u_1) + ... + f_T(u_\Delta) + [\mu(m')
- \mu(m)]$. Weights are degree-monotone in $\langle\mu,d\rangle$, so
we have $\mu(m') - \mu(m) \ge 0$, and, since $\Delta> 0$, from Note
\ref{note_positive} it follows that $b > 0$.

By construction of $T'$ we have $f_{T'}(m') = f_T(m)$, $f_{T}(m') - f_{T'}(m)= b$. Therefore, by Lemma
\ref{lemma_Zhang_b},
$$\mathbf{y} := (f_{T'}(m)) =(f_{T}(m')-b)\succ^w (f_T(m')) =: \mathbf{x}.$$
Weights of groups subordinated to all other vertices do not change,
so, by analogy with Lemmas \ref{lemma_major1}-\ref{lemma_major5a},
by Lemma \ref{lemma_Zhang_xy} we obtain
$\mathbf{f}(T')\succ^w\mathbf{f}(T)$.
\end{proof}
\end{lemma}

Please note that only Lemmas \ref{lemma_major4} and \ref{lemma_major6} require degree-monotonicity of the generating
tuple $\langle\mu,d\rangle$.

As we show below, conditions of Lemmas
\ref{lemma_major1}-\ref{lemma_major6} are never satisfied for
directed Huffman trees (an only for directed Huffman trees), and the
above transformations cannot result in a tree with the vector of
subordinate groups' weights majorizing the one of some directed
Huffman tree.

\subsection{Huffman Trees and Majorization}

Let us define the following useful tree transformations.

\begin{definition}
A directed tree $T \in \mathcal{R}(\mu, d)$ \emph{induces} the Huffman tree $H$, if $H \in \mathcal{RH}(\mu, d)$. A
directed tree $H' \in \mathcal{R}(\mu, d)$ is an \emph{augmentation} of a Huffman tree for the $m$-rollup of the tree
$T \in \mathcal{R}(\mu, d)$ if $H'$ is obtained by joining, firstly, the Huffman tree $\underline{H}$ induced by an
$m$-rollup $\underline{T}$ of the directed tree $T$, and, secondly, the contraction $R \in \mathcal{R}$ of the tree $T$
to the vertex set $V(R) := g_T(m)$, i.e. $E(R):=E(T)\cap (g_T(m)\times g_T(m)), V(H') := V, E(H') := E(\underline{H})
\cup E(R), \mu_{H'}(v) := \mu(v), v \in V$.
\end{definition}

\begin{note}\label{note_augm}
If the vertex $m\in M$ is not the root of $T$, then, clearly, $\mathbf{f}(H') = (\mathbf{f}(R), f_T(m),
\mathbf{f}(\underline{H}))_\uparrow$. Moreover, formula (\ref{eq_weights_monotone}) from Lemma
\ref{lemma_Huffman_monotonicity} holds for those vertices of the augmented tree $H'$, which also belong to the Huffman
tree $\underline{H}$ induced by the $m$-rollup of $T$.
\end{note}

\begin{theorem}\label{theorem_major}
If weights are degree-monotone in a generating tuple $\langle\mu,
d\rangle$, and $H\in \mathcal{RH}(\mu, d)$ is a directed Huffman
tree, then for any directed tree $T \in \mathcal{R}(\mu, d)$
$\mathbf{f}(H)\succeq^w \mathbf{f}(T)$.
\begin{proof}
Let us employ induction on the number of internal vertices $q$. If $q = 1$, the statement of the theorem is
straightforward, since the collection $\mathcal{R}(\mu, d)$ consists of the sole tree (a directed star). Assume the
theorem is valid for all directed trees with less than $q$ internal vertices. Let us prove that it is also valid for
directed trees with $q$ internal vertices.

The relation $\succ^w$ is a strict partial ordering on the set of
$(q - 1)$-dimensional vectors, and, hence, a strict partial ordering
on a narrower set of vectors of subordinate groups' weights of all
directed trees from $\mathcal{R}(\mu, d)$. Therefore, the set
$$\bar{\mathcal{R}}(\mu, d) := \{T \in \mathcal{R}(\mu, d): \nexists T' \in \mathcal{R}(\mu, d)\textrm{ such that }\mathbf{f}(T') \succ^w \mathbf{f}(T)\}$$
of trees whose vector of subordinate groups' weights is ``maximal''
with respect to the partial ordering $\succ^w$, is not empty.
Without loss of generality suppose that $T \in
\bar{\mathcal{R}}(\mu, d)$.

Among all lower stars of the tree $T$, one or more has the least total weight. One or more centers of these
least-weight lower stars has the least degree. Let $\underline{v} \in M$ be one of these least-degree centers having
the least vertex weight $\mu(\cdot)$, and let $\underline{S}$ be the $\underline{v}$-centered lower star in $T$.

Note that, since $q> 1$, the vertex $\underline{v}$ (being a center of a lower star) cannot be the root of $T$. The
following four steps prove that the star $\underline{S}$ is minimal for the tuple $\langle\mu, d\rangle$, i.e., that
$f_T(\underline{v}) = \underline{f}(\mu, d)$. Below the shorthand notation $d := \min_{u\in M} d(u)$ is used.

\textbf{Step I.} First we prove that the tree $T$ contains a lower star with $d - 1$ pendent vertices. Suppose, by
contradiction, that centers of all lower stars in $T$ have more than $d-1$ pendent vertices and, thus, have degree
greater than $d$. Then the tree $T$ must contain a vertex $m' \in M$ of degree $d_T(m') = d$, which has an inbound arc
from some vertex $m \in M$ of degree $d_T(m)
> d$. But Lemma \ref{lemma_major6} says that then the tree $T' \in \mathcal{R}(\mu, d)$ exists, such that $\mathbf{f}(T') \succ^w \mathbf{f}(T)$,
so $T$ cannot belong to the collection $\bar{\mathcal{R}}(\mu, d)$.
The obtained contradiction proves that $T$ contains some lower star
(denote it with $S$) containing $d - 1$ pendent vertices. Let $v\in
M$ be the center of $S$.

\textbf{Step II.} Let us prove that the star $\underline{S}$ has exactly $d - 1$ leaves. Suppose, by contradiction,
that $|L(\underline{S})|> d - 1$. In particular, this implies that $S \neq \underline{S}$ and $f_T(v) >
f_T(\underline{v})$ ($f_T(v) \ge f_T(\underline{v})$ by construction of the star $\underline{S}$, and the case of
$f_T(v) = f_T(\underline{v})$ contradicts the fact that $\underline{v}$ has the least degree among all least-weight
lower stars).

It is clear that the star $S$ is still a lower star in a
$\underline{v}$-rollup of the tree $T$, so, let $H^\clubsuit \in
\mathcal{R}(\mu, d)$ stand for the augmentation of a Huffman tree
$\uwave{H}$ induced by the $(v, \underline{v})$-rollup $\uwave{T}$
of the tree $T$. Since neither $v$, nor $\underline{v}$, are the
roots of $T$, by Definition \ref{def_rollup} we have $\mathbf{f}(T)
= (f_T(v), f_T(\underline{v}), \mathbf{f}(\uwave{T}))$. By Note
\ref{note_augm}, $\mathbf{f}(H^\clubsuit) = (f_T(v),
f_T(\underline{v}), \mathbf{f}(\uwave{H}))$. By inductive assumption
we have $\mathbf{f}(\uwave{H}) \succeq^w \mathbf{f}(\uwave{T})$, so,
by Lemma \ref{lemma_Zhang_xy}, $\mathbf{f}(H^\clubsuit) \succeq^w
\mathbf{f}(T)$. Since $T \in \bar{\mathcal{R}}(\mu, d)$, the case of
$f(H^\clubsuit) \succ^w f(T)$ is impossible, and, thus,
$\mathbf{f}(H^\clubsuit) = \mathbf{f}(T)$, i.e., $H^\clubsuit \in
\bar{\mathcal{R}}(\mu, d)$.

Definitely, disjoint paths $(v, m_1, ..., m_k, m)$ and
$(\underline{v}, m_1', ..., m_l', m)$ to some vertex $m\in M$
present in $H^\clubsuit$, where $k, l \ge 0$. Again recall Note
\ref{note_augm}: since $f_T(v)> f_T(\underline{v})$, formula
(\ref{eq_weights_monotone}) makes $f_{H^\clubsuit} (m_i) >
f_{H^\clubsuit} (m_i'), i = 1, ..., \min[k, l]$. It also follows
from (\ref{eq_weights_monotone}) that $k \le l$, since otherwise
$f_{H^\clubsuit}(m_{l+1}) > f_{H^\clubsuit}(m)$, which is
impossible, as $m_{l+1} \in g_{H^\clubsuit}(m)$. Thus, the
considered pair of paths satisfies conditions of Lemma
\ref{lemma_major4}, and a tree exists whose vector of subordinate
groups' weights majorizes the appropriate vector of the tree
$H^\clubsuit$, which contradicts the fact that $H^\clubsuit \in
\bar{\mathcal{R}}(\mu, d)$.

The obtained contradiction proves that the star $\underline{S}$ has $d - 1$ pendent vertices.

\textbf{Step III.} Let us prove that the vertex $\underline{v}$ (the center of the star $\underline{S}$) has the least
weight $\mu(\cdot)$ in the set $M$. Assume, by contradiction, that a vertex $\underline{\underline{v}}\in M$ exists
such that $\mu(\underline{\underline{v}}) < \mu(\underline{v})$. Since weights are degree-monotone in the tuple
$\langle\mu,d\rangle$, this implies that $d(\underline{v}) = d(\underline{\underline{v}})=d$. By construction of the
vertex $\underline{v}$ we have $f_T(\underline{v}) \le f_T(\underline{\underline{v}})$. Moreover, we can discard the
case of $f_T(\underline{v}) = f_T(\underline{\underline{v}})$, since then the vertex $\underline{\underline{v}}$ would
be the center of a lower star, and, because $d(\underline{v}) = d(\underline{\underline{v}})$, we would not have
$\mu(\underline{\underline{v}}) < \mu(\underline{v})$ by construction of the vertex $\underline{v}$. Consequently, only
the case of $f_T(\underline{\underline{v}})> f_T(\underline{v})$ is left.

If $\underline{v}\in g_T(\underline{\underline{v}})$, a path exists
from $\underline{v}$ to $\underline{\underline{v}}$ in $T$, and, by
Lemma \ref{lemma_major5a}, $T\notin\bar{\mathcal{R}}(\mu, d)$.
Otherwise consider an augmentation $H^\diamondsuit \in
\mathcal{R}(\mu, d)$ of a Huffman tree $\underline{\underline{H}}$
induced by the $(\underline{\underline{v}}, \underline{v})$-rollup
$\underline{\underline{T}}$ of the tree $T$. If
$\underline{\underline{R}}$ is a contraction of $T$ to the vertex
set $g_T(\underline{\underline{v}})$, then, by Note \ref{note_augm},
we have $\mathbf{\mathbf{f}}(H^\diamondsuit) = (f_T(\underline{v}),
\mathbf{f}(\underline{\underline{R}}),
f_T(\underline{\underline{v}}),
\mathbf{f}(\underline{\underline{H}}))_\uparrow$. On the other hand,
by Definition \ref{def_rollup}, $\mathbf{f}(T) =
(f_T(\underline{v}), \mathbf{f}(\underline{\underline{R}}),
f_T(\underline{\underline{v}}),
f(\underline{\underline{T}}))_\uparrow$. By inductive assumption,
$\mathbf{f}(\underline{\underline{H}}) \succeq^w
\mathbf{f}(\underline{\underline{T}})$, i.e., by Lemma
\ref{lemma_Zhang_xy}, $\mathbf{f}(H^\diamondsuit) \succeq^w
\mathbf{f}(T)$. Since, by assumption, $T \in \bar{\mathcal{R}}(\mu,
d)$, the case of $\mathbf{f}(H^\diamondsuit) \succ^w \mathbf{f}(T)$
is impossible, so, $\mathbf{f}(H^\diamondsuit) = \mathbf{f}(T)$, and
$H^\diamondsuit \in \bar{\mathcal{R}}(\mu, d)$. Again,
$H^\diamondsuit$ contains disjoint paths $(v, m_1, ..., m_k, m)$ and
$(v, m_1', ..., m_l', m)$, $k, l \ge 0$ to some vertex $m \in M$.
Since $f_T(\underline{\underline{v}})> f_T(\underline{v})$, applying
formula (\ref{eq_weights_monotone}) we deduce that $k \le l$,
$f_{H^\diamondsuit}(m_i)
> f_{H^\diamondsuit}(m_i')$, $i = 1, ..., k$. Then, by Lemma \ref{lemma_major5}, the vector of subordinate group weights
of $H^\diamondsuit$ is majorized by the appropriate vector of some
tree from $\mathcal{R}(\mu,d)$, and $H^\diamondsuit$ cannot be in
$\bar{\mathcal{R}}(\mu, d)$. The obtained contradiction proves that
the vertex $\underline{v}$ has the least weight in $M$.

\textbf{Step IV.} Now to prove the minimality of the star $\underline{S}$ we are left to show that its pendent vertices
have the least possible weights $\mu(\cdot)$. Assume, by contradiction, that it is not true, i.e., such vertices $w \in
W\backslash L(\underline{S})$ and $w' \in L(\underline{S})$ exist that $\mu(w) < \mu(w')$. The vertex $w$ has an
outbound arc in $T$ to some vertex $\tilde{v}\in M$. There are two possible alternatives:
\begin{enumerate}
    \item $\underline{v} \in g_T(\tilde{v})$. By assumption, $f_T(w) = \mu(w) < \mu(w') = f_T(w')$ and thus, by
    Lemma \ref{lemma_major3} we conclude that $T \notin \bar{\mathcal{R}}(\mu, d)$, which contradicts the above assumption.
    \item $\underline{v} \notin g_T(\tilde{v})$. Let $H^\heartsuit$ be the augmentation of a Huffman tree induced
    by a $(\underline{v}, \tilde{v})$-rollup of the tree $T$.
    By analogy to the step II we show that $H^\heartsuit \in \bar{\mathcal{R}}(\mu, d)$.

By construction, there are disjoint paths $(w, m_1, ..., m_k, m)$
and $(w', m_1', m_2', ..., m_l', m)$ (where $k \ge 1$, $l \ge 1$,
$m_1 = \tilde{v}, m_1' = \underline{v}$), in $H^\heartsuit$ to some
vertex $m \in M$. We have $f_T(\underline{v}) \le f_T(\tilde{v})$ by
construction of the vertex $\underline{v}$. If this inequality is
strict, then we also have $f_{H^\heartsuit}(\underline{v}) <
f_{H^\heartsuit}(\tilde{v})$ and, using formula
(\ref{eq_weights_monotone}), conclude that $k \le l$,
$f_{H^\heartsuit}(m_i)> f_{H^\heartsuit}(m_i')$, $i = 1, ..., k$.
Since $\mu(w) < \mu(w')$, Lemma \ref{lemma_major1} predicates the
existence of a tree, whose vector of subordinate groups' weights
majorizes the appropriate vector of the tree $H^\heartsuit$, which
contradicts to the assumption that $H^\heartsuit \in
\bar{\mathcal{R}}(\mu, d)$.

In case of $f_T(\underline{v}) = f_T(\tilde{v})$ we cannot use formula (\ref{eq_weights_monotone}) to compare
subordinate groups' weights of elements of both chains, since all possible alternatives of $k = 1$, or $l = 1$, or any
sign of the expression $f_{H^\heartsuit}(m_2) - f_{H^\heartsuit}(m_2')$ in case of $k, l \ge 2$ are possible. On the
other hand, if $f_{H^\heartsuit}(m_2)> f_{H^\heartsuit}(m_2')$, then formula (\ref{eq_weights_monotone}) can be used to
show that $k \le l$, $f_{H^\heartsuit}(m_i) > f_{H^\heartsuit}(m_i')$, $i = 2, ..., k$. In case of the opposite
inequality, $f_{H^\heartsuit}(m_2) < f_{H^\heartsuit}(m_2')$, formula (\ref{eq_weights_monotone}) says that, by
contrast, $k \ge l$, $f_{H^\heartsuit}(m_i) < f_{H^\heartsuit}(m_i')$, $i = 2, ..., l$. Repeating this argument through
the chain, we see that only two alternatives are possible:
\begin{itemize}
    \item $1\le p \le k \le l$, $f_{H^\heartsuit}(m_i) = f_{H^\heartsuit}(m_i')$, $i = 1, ..., p$,
    $f_{H^\heartsuit}(m_i) > f_{H^\heartsuit}(m_i')$, $i = p + 1, ..., k$. Since $\mu(w) < \mu(w')$,
    Lemma \ref{lemma_major1} gives $H^\heartsuit \notin \bar{\mathcal{R}}(\mu, d)$, which is a contradiction.
    \item $1\le p \le l \le k$, $f_{H^\heartsuit}(m_i) = f_{H^\heartsuit}(m_i')$, $i = 1, ..., p$,
    $f_{H^\heartsuit}(m_i) < f_{H^\heartsuit}(m_i')$, $i = p + 1, ..., l$. In this case the same conclusion that
    $H^\heartsuit \notin \bar{\mathcal{R}}(\mu, d)$ follows from Lemma \ref{lemma_major2}.
\end{itemize}
\end{enumerate}

The obtained contradictions prove the minimality of the star $\underline{S}$ having the center $\underline{v}$. In
other words, we have $f_T(\underline{v}) = \underline{f}(\mu, d)$.

Let $H^\spadesuit$ be the augmentation of the Huffman tree
$\underline{H}$ induced by the $\underline{v}$-rollup
$\underline{T}$ of the tree $T$. Since $|M(\underline{H})| =
|M(\underline{T})| = q - 1$, by inductive assumption we have
$f(\underline{H}) \succeq^w f(\underline{T})$. Since $\mathbf{f}(T)
= (\underline{f}(\mu, d), \mathbf{f}(\underline{T}))$,
$f(H^\spadesuit) = (\underline{f}(\mu, d),
\mathbf{f}(\underline{H}))$, by Lemma \ref{lemma_Zhang_xy} obtain
$\mathbf{f}(H^\spadesuit) \succeq^w \mathbf{f}(T)$. As
$H^\spadesuit$ is constructed by adding a minimal star
$\underline{S}$ to the Huffman tree $\underline{H}$, by Lemmas
\ref{lemma_star_seq_weight} and \ref{lemma_same_weights_is_Huffman}
it appears to be a Huffman tree itself, i.e., $H^\spadesuit\in
\mathcal{RH}(\mu, d)$. Then Lemma
\ref{lemma_all_Huffman_same_weights} says that
$\mathbf{f}(H^\spadesuit) = \mathbf{f}(H)$, and, therefore,
$\mathbf{f}(H) \succeq^w \mathbf{f}(T)$.
\end{proof}
\end{theorem}

\section{Huffman Tree Minimizes Wiener Index}\label{section_minimizing_index}

\begin{definition}A directed tree $T \in \mathcal{R}(\mu, d)$ with the vector $(f_1, ..., f_{q-1})$
of subordinate groups' weights is called a \emph{proper tree} if $f_i \le \bar{\mu}/2$, $i = 1, ..., q - 1$.
\end{definition}

\begin{lemma}\label{lemma_proper_tree_exists}
 Each tree $T \in \mathcal{T}(\mu, d)$ has a corresponding proper tree, and vice versa.
\begin{proof}
For a tree with $q < 2$ internal vertices the lemma is trivial, since the vector of subordinate groups' weights is
empty, therefore, suppose that $q \ge 2$. Consider a vertex $u \in V(T)$ with incident edges $uv_1, ..., uv_d \in
E(T)$, where $d := d_T(u)$. Deleting the vertex $u$ and the edges $uv_1, ..., uv_d$ we break the tree $T$ into disjoint
components  $T_1, ..., T_d$.

Let us prove that in any tree $T \in \mathcal{T}(\mu, d)$ deletion of some vertex $v\in V(T)$ results in components of
the weight no more than $\bar{\mu}/2$. Assume, by contradiction, that for every vertex $u \in V(T)$ such an incident
vertex $\pi(u) \in V(T)$ exists that deletion of the edge $u\pi(u)$ gives rise to the component including the vertex
$\pi(u)$ and having the weight greater than $\bar{\mu}/2$. Clearly, the weight of the second component (the one
including the vertex $u$) does not exceed $\bar{\mu}/2$. Then, since the tree $T$ has finite number of vertices, it
inevitably contains a cycle $(u_1, ..., u_k, u_1)$, where $k>2$, $u_{i+1} = \pi(u_i)$, $i = 1, ..., k - 1$, $u_1 =
\pi(u_k)$, which contradicts to the fact that $T$ is a tree.

If $v$ is an internal vertex, then we choose the root $r = v$, otherwise let the vertex incident to $v$ be the root $r$
of the corresponding directed tree $P\in \mathcal{R}(\mu,d)$ ($r \in M(T)$, since $q \ge 1$). One can easily see that
$P$ is a proper tree. The inverse statement is trivial.
\end{proof}
\end{lemma}

Let $\mathcal{P} \subseteq \mathcal{R}$ stand for the collection of all proper directed trees, and let
$\mathcal{P}(\mu, d) \subseteq \mathcal{R}(\mu, d)$ be the collection of proper trees, which correspond to the trees
from $\mathcal{T}(\mu, d)$.

\begin{lemma}\label{lemma_Huffman_proper}
If weights are degree-monotone in the tuple $\langle\mu, d\rangle$, all directed Huffman trees from $\mathcal{RH}(\mu,
d)$ are proper trees.
\begin{proof}
Consider a directed Huffman tree $H \in \mathcal{RH}(\mu, d)$ with a star sequence $S_1, ..., S_q$ and let $m_1, ...,
m_q$ be the centers of stars $S_1, ..., S_q$ respectively. To prove the lemma we are enough to show that $f_{q-1}(H) =
f_H(m_{q-1}) \le \bar{\mu}/2$. Assume, by contradiction, that $f_H(m_{q-1}) > \bar{\mu}/2$. Since $f_H(m_q) =
\bar{\mu}$, we have
$$\sum_{v\in L(S_q)\backslash\{m_{q-1}\}}f_H(v)=\bar{\mu}-f_H(m_{q-1})-\mu(m_q)<$$
$$<\bar{\mu}/2-\mu(m_q)<f_H(m_{q-1})-\mu(m_q)=\mu(m_{q-1})+\sum_{v\in L(S_{q-1})}f_H(v)-\mu(m_q).$$
Since $H \in \mathcal{RH}(\mu, d)$, by construction of the Huffman tree we have $\mu(m_{q-1}) \le \mu(m_q)$, so,
\begin{equation}\label{eq_non_monotone}
\sum_{v\in L(S_q)\backslash\{m_{q-1}\}}f_H(v)<\sum_{v\in L(S_{q-1})}f_H(v).
\end{equation}

From the fact that $H \in \mathcal{RH}(\mu, d)$ and from degree-monotonicity of weights in $\langle\mu, d\rangle$ it
follows that $d_H(m_{q-1}) \le d_H(m_q)$. Introduce the notation $\Delta := d_H(m_{q-1}) - 1$. Choose any $\Delta$
vertices from the set $L(S_q) \backslash \{m_{q-1}\}$ and denote them with $v_1, ..., v_\Delta$. Transform the tree $H$
by redirecting outbound arcs of vertices $v_1, ..., v_\Delta$ to the vertex $m_{q-1}$ and by redirecting all $\Delta$
outbound arcs from the vertices of the set $L(S_{q-1})$ to the vertex $m_q$. Vertex degrees and weights do not change
during the transformation, thus, the transformed tree $H'$ belongs to $\mathcal{R}(\mu, d)$. Using inequality
(\ref{eq_non_monotone}) it is easy to show that the weight of the group subordinated to the vertex $m_{q-1}$ decreases
by
$$b:= \sum_{v\in L(S_{q-1})}f_H(v)-\sum_{i=1}^\Delta f_H(v_i) \ge \sum_{v\in L(S_{q-1})} f_H(v) - \sum_{v\in L(S_q)\backslash\{m_{q-1}\}}f_H(v)>0,$$
while the weights of groups subordinated to all other vertices in a
tree do not change. Then, by Lemma \ref{lemma_Zhang_b}, $f(H')
\succ^w f(H)$, which is impossible by Theorem \ref{theorem_major}.
The obtained contradiction proves that $f_H(m_{q-1}) \le
\bar{\mu}/2$ and $H$ is a proper tree.
\end{proof}
\end{lemma}

\begin{theorem}\label{theorem_concave}
Let a function $\chi(x)$ be concave and increasing on the range $x
\in [0, \bar{\mu}/2]$. If $H\in\mathcal{RH}(\mu, d)$, $T \in
\mathcal{P}(\mu, d)$, then $\sum_{v\in V\backslash
\{r\}}\chi(f_H(v))\le\sum_{v\in V\backslash \{r'\}}\chi(f_T(v))$,
where $r$ and $r'$ are the roots of the trees $H$ and $T$
respectively. The equality is possible only if $T\in
\mathcal{RH}(\mu, d)$.
\begin{proof}
By Theorem \ref{theorem_major}, if $H\in \mathcal{RH}(\mu, d)$, $T
\in \mathcal{P}(\mu, d)\subseteq \mathcal{R}(\mu, d)$, then
$\mathbf{f} := \mathbf{f}(H) \succeq^w \mathbf{f'} :=
\mathbf{f}(T)$. From Lemma \ref{lemma_Huffman_proper} we know that
$H\in\mathcal{P}(\mu, d)$, so both $f_i$ and $f_i'$ ($i = 1, ..., q
- 1$), belong to the range $[0,\bar{\mu}/2]$, where the function
$\chi(x)$ is increasing and concave. Then, by Lemma
\ref{lemma_Zhang_concave}, $\sum_{i=1}^{q-1}\chi(f_{[i]}) \le
\sum_{i=1}^{q-1}\chi(f_{[i]}')$, and the equality is possible only
if $\mathbf{f} = \mathbf{f'}$. Since trees from $\mathcal{R}(\mu,d)$
differ only in weights of groups subordinated to internal vertices,
we immediately obtain the desired inequality. In case of equality,
Lemma \ref{lemma_same_weights_is_Huffman} proves
 that $T \in \mathcal{RH}(\mu, d)$.
\end{proof}
\end{theorem}

Now we are ready to prove the \textbf{Theorem \ref{theorem_main}}.
\begin{proof}
From Lemma \ref{lemma_proper_tree_exists}, each tree in $T \in
\mathcal{T}(\mu, d)$ has a corresponding proper directed tree $P$,
and vice versa. From equation (\ref{eq_VWWI_directed}) it follows
that trees $P$ and $T$ share the same value of the Wiener index, so,
if $\mathcal{P}^* (\mu,d):=\textrm{Argmin}_{T\in \mathcal{P}(\mu,d)}
VWWI(T)$, then the collection $\mathcal{T}^*(\mu,d)$ of
vertex-weighted trees minimizing the Wiener index is a collection of
trees corresponding to trees from $\mathcal{P}^*(\mu,d)$. Since
\emph{VWWI} satisfies the conditions of Theorem
\ref{theorem_concave} with function $\chi(x)=x(\bar{\mu}-x)$, which
is concave and increasing on the range $[0, \bar{\mu}/2]$, we deduce
that $\mathcal{P}^*(\mu,d) = \mathcal{RH}(\mu, d)$ and, since
$\mathcal{TH}(\mu, d)$ is the collection of all trees corresponding
to directed trees from $\mathcal{RH}(\mu, d)$, we finally deduce
that $\mathcal{T}^*(\mu,d) = \mathcal{TH}(\mu, d)$.
\end{proof}

\section{Conclusion}

In the conclusion, let us discuss possible applications and extensions of the considered model. In
\cite{Goubko14,GoubkoGutman14,GoubkoReti14} a technique was suggested to optimize abstract degree-based topological
indices of the form $C_{deg}(G):=C_1(G)+C_2(G)$, where
$$
C_1(G)=\sum_{v\in V(G)}c_1(d_G(v)), \qquad C_2(G)=\sum_{uv\in E(G)}c_2(d_G(u), d_G(v)),
$$
over the set of trees with the given number of pendent vertices. Together with the results of this paper the technique
from \cite{Goubko14} can be seen as a step towards optimization of joint linear combinations of degree-based and
distance-based indices over the set of trees with the given total number of vertices or the given number of pendent
vertices.

For the fixed degree sequence optimization of the linear combination of $C_1(G)$ and $VWWI(G)$ reduces to building a
Huffman tree, and then we are just to find an optimal degree sequence, which is an integer program with linear
constraint (\ref{tree_degrees_identity}). The efficient algorithms for joint index optimization problems would
contribute to the methods of designing materials with extremal characteristics.

In the definition of $VWWI(G)$ each distance from $u\in V(G)$ to $v\in V(G)$ is weighted with the product
$\mu_G(u)\mu_G(v)$ of positive vertex weights. The obvious extension assumes endowing each path in a tree (i.e., each
pair $u, v \in V$) with its own weight $\mu_G(u,v)$. Such an extended index would give more freedom to build models
relating physical and chemical properties of substances to the topology of their molecules. For instance, we would be
able to assign independently unique weights to $\textrm{OH-OH}$, $\textrm{C-OH}$ and $\textrm{C-C}$ paths in polyhydric
alcohol molecules.

This settings seems to be closer to NP-hard problems of hierarchical graph clustering (see \cite{Shaefer07} for a
survey) and merits a more detailed study in the future.

\end{document}